\documentclass[12pt]{amsart} \usepackage{amssymb} 
\newtheorem{thm}{Theorem}[section]
\newtheorem{cor}[thm]{Corollary} \newtheorem{lem}[thm]{Lemma}
\newtheorem{conj}[thm]{Conjecture}

\newtheorem{prop}[thm]{Proposition}

\newtheorem{definition}[thm]{Definition}

\newtheorem{example}[thm]{Example}
\newenvironment{exmp}{\begin{example}\rm}{\end{example}}
\newtheorem{remark}[thm]{Remark}
\newenvironment{rem}{\begin{remark}\rm}{\end{remark}}

\font\tenmsb=msbm10 \font\sevenmsb=msbm7 \font\fivemsb=msbm5

\newfam\msbfam \textfont\msbfam=\tenmsb
\scriptfont\msbfam=\sevenmsb \scriptscriptfont\msbfam=\fivemsb

\font\teneufm=eufm10 \font\seveneufm=eufm7 \font\fiveeufm=eufm5
\newfam\eufmfam \textfont\eufmfam=\teneufm
\scriptfont\eufmfam=\seveneufm
\scriptscriptfont\eufmfam=\fiveeufm

\newcommand{\zwei}[2]{\left[
    \begin{array}{c} #1 \\ #2
    \end{array} \right]} \newcommand{\vier}[4]{\left[
    \begin{array}{ccc} #1 &\;& #2 \\ #3 &\;& #4 \end{array}
  \right]}

\newcommand{\F}{{\mathbb F}} \newcommand{\C}{{\mathbb C}}
\newcommand{\R}{{\mathbb R}}

\title[Some Remarks on Real and Complex Output Feedback]{Some Remarks on\\ Real and Complex Output Feedback}
\author{Joachim Rosenthal}
\address{\hskip-\parindent
Joachim Rosenthal\\
Department of Mathematics\\
University of Notre Dame\\ 
Notre  Dame, IN 46556, USA}
\email{Rosenthal.1@nd.edu}
\author{Frank Sottile}
\address{\hskip-\parindent
Frank Sottile\\
Department of Mathematics\\
        University of Toronto\\
        100 St.~George Street\\
	Toronto, Ontario  M5S 3G3\\
	CANADA} 
\curraddr{MSRI\\
	1000 Centennial Drive\\
	Berkeley CA, 94720\\
	USA}
\email{sottile@msri.org, sottile@math.toronto.edu}

\thanks{First author supported in part by NSF grant
DMS-94-00965.} 
\thanks{Second author supported in part by NSERC grant  OGP0170279 \ 
and NSF grant DMS-9022140}
\date{January 28, 1997}

\begin{document} 
\begin{abstract}
   We provide some new necessary and
  sufficient conditions which guarantee arbitrary pole placement
  of a particular linear system over the complex numbers.
  We exhibit a non-trivial real
 linear system which is not controllable by real static output
 feedback and discuss a conjecture from algebraic geometry
 concerning the existence of real linear systems for which all
 static feedback laws are real.
\end{abstract}

\maketitle

\section{Preliminaries}

Let $\F$ be an arbitrary field and let $m,p,n$ be fixed positive
integers. Let $A,B,C$ be matrices with entries in $\F$ of sizes
$n\times n$, $n\times m$, and $p\times n$ respectively. Identify
the space of monic polynomials having degree $n$,
$$
s^n+a_{n-1}s^{n-1}+\cdots+a_1s+a_0\in \F[s],
$$
with the vector space $\F^n$. In its simplest form, the static
output pole placement problem asks for conditions on the matrices
$A,B,C$ which guarantee that the pole placement map
\begin{equation}
  \label{pole-map}
  \chi_{(A,B,C)}:\, \F^{mp}\longrightarrow\F^n,\hspace{3mm}
  K\longmapsto \det(sI-A-BKC)
\end{equation}
is surjective. A dimension argument shows the necessity of
$mp\geq n$.  This question has been studied intensively and we
refer to the survey articles~\cite{by89,ro97} and the recent
papers~\cite{ki95p,le95,ro95,ro96,wa92,wa96} for details.  We
summarize some of the most important results.

A matrix pair $(A,B)$ defined over a field $\F$ is {\em
  controllable} if the matrix pencil $\left[ sI-A \mid B\right]$
is left coprime. Equivalently, if the full size minors of the pencil
$\left[ sI-A \mid B\right]$ have no common non-trivial polynomial
factor.  Similarly, a matrix pair $(A,C)$ is {\em observable} if
the matrix pencil $\zwei{sI-A}{C}$ is right coprime. Then we
have:
\begin{lem}
  $\chi_{(A,B,C)}$ is surjective only if $(A,B)$ is a
  controllable pair and $(A,C)$ is an observable pair.
\end{lem}
\begin{proof}
  The following identity immediately establishes the claim:
  $$
  \det(sI-A-BKC)=\det\vier{sI-A}{-B}{-KC}{I}=
  \det\vier{sI-A}{-BK}{-C}{I}
  $$
\end{proof}

The necessary conditions $mp\geq n$, controllability, and
observability are not sufficient to guarantee arbitrary pole
assignability.  When $p=1$, the following straightforward lemma
provides exact conditions for arbitrary pole assignability over
any field $\F$.
\begin{lem}\label{p=1suff}
  Let $p=1$ and let $d^{-1}(s)(n_1(s),\ldots,n_m(s))$ be a left
  coprime factorization of the transfer function $C(sI-A)^{-1}B$.
  Then the pole placement map (\ref{pole-map}) is surjective if
  and only if $n_1(s),\ldots,n_m(s)$ span the vector space of
  polynomials of degree at most $n-1$.
\end{lem}
One readily establishes a similar result when $m=1$.
Lemma~\ref{p=1suff} gives algebraic conditions on the set of
systems parameters. To make this precise, identify the set of
matrices $(A,B,C)$ having fixed sizes $n\times n$, $n\times m$,
and $p\times n$ with the vector space $V:=\F^{n(n+m+p)}$. Recall
that a subset $G\subset V$ is {\em generic} if a non-trivial
polynomial vanishes on its complement $V\setminus G$.  Thus
Lemma~\ref{p=1suff} implies that if $p=1$ and $m\geq n$, then the
set of systems which can be arbitrarily pole assigned forms a
generic set.

Since non-controllable systems $(A,B,C)$ cannot be arbitrarily
pole assigned, pole placement results are often restricted to a
generic class of systems. If the base field $\F$ is the real
numbers $\R$ or the complex numbers $\C$, then a generic set
$G\subset V$ is open and dense with respect to the usual
Euclidean topology, and its complement $V\setminus G$ has measure
zero.

If the pole placement map $\chi$ is surjective for a generic set
of systems and some fixed base field $\F$ we will say in short
that $\chi$ is {\em generically surjective}.

The major results are as follows:
\begin{thm}[Brockett and Byrnes~\cite{br81}]   \label{BrBy}
  If the base field $\F$ is algebraically closed and if $mp\geq
  n$ then $\chi$ is generically surjective. Moreover if $mp=n$
  then for a generic set of systems the cardinality of
  $\chi^{-1}(\phi)$ (when counted with multiplicity) is
  independent of the closed loop polynomial $\phi\in \F^n$ and is
  equal to
\begin{equation}                       \label{1}
d(m,p)= \frac{1!2!\cdots (p-1)!(mp)!}{m!(m+1)!\cdots(m+p-1)!} 
\end{equation}
\end{thm}

Since $mp\geq n$ is necessary for $\chi$ to be surjective,
Theorem~\ref{BrBy} gives the best possible bound when the base
field $\F$ is algebraically closed.

The number $d(m,p)$ is the degree of the Grassmann variety, which
was computed in the last century by Schubert~\cite{sh91}.
Although the real numbers $\R$ are not algebraically closed and
Theorem~\ref{BrBy} therefore does not apply one still has the
following Corollary:
\begin{cor}\label{oddR}
  If\/ $\F=\R$, $mp=n$, and $d(m,p)$ is odd, then $\chi$ is
  generically surjective.
\end{cor}
\begin{proof}
  If $(A,B,C)$ are real matrices then the set $\chi^{-1}(\phi)$
  is closed under complex conjugation for every closed loop
  polynomial $\phi\in \R^n$. Therefore, for a generic set of systems,
  $\chi^{-1}(\phi)$ contains a real point for each $\phi$.
\end{proof}

As an example, consider the case $\F=\R$, $m=2$, $p=3$ and $n=6$.
Here, $d(2,3)=5$.  At least one of the 5 points $\chi^{-1}(\phi)$
is real, so $\chi$ is generically surjective even over the reals.

Berstein determined when $d(m,p)$ is odd.
\begin{prop}[Berstein~\cite{be76}]
  The number $d(m,p)$ is odd if and only if $\min (m,p)=1$ or
  $\min (m,p)=2$ and $\max (m,p)=2^t-1,$ where $t$ is a positive
  integer.
\end{prop}
When $d(m,p)$ is even, the best known sufficiency result over the
reals is due to Wang:
\begin{thm}[Wang~\cite{wa92}]
  If\/ $\F=\R$ and $mp > n$, then $\chi$ is generically
  surjective.
\end{thm}

For an elementary direct proof of this important sufficiency
result we refer to~\cite{ro95}.

For generic surjectivity over the reals, there is a difference of
one degree of freedom between sufficiency ($mp>n$) and necessity
($mp\geq n$).  As we already noted, $mp\geq n$ is sufficient if
$d(m,p)$ is an odd number.  One may ask if $mp\geq n$ might be
always sufficient?
\begin{prop}[Willems and Hesselink~\cite{wi78}]
  If\/ $\F=\R$ and if $m=p=2$ and $n=4$ then there is an open
  Euclidean neighborhood $U\subset V=\R^{32}$ having the property
  that $\chi_{(A,B,C)}$ is not surjective if $(A,B,C)\in U$. In
  particular $\chi$ is not generically surjective.
\end{prop}

It has been conjectured by S.-W. Kim that $m=p=2$, $n=4$ is the
only case where $mp=n$ is not a sufficient condition for $\chi$
to be generically surjective over the reals. In the next section
we exhibit a counterexample.

\section{Main Results}

The result by Brockett and Byrnes provides a sufficiency result
for a generic set of systems.  We provide exact conditions which
guarantee that a particular plant $(\bar{A},\bar{B},\bar{C})$ is
arbitrarily pole assignable.  Our approach is geometric,
utilizing the central projection of the Grassmann variety induced
by the pole placement map~\cite{br81,wa96}.

Let $D^{-1}(s)N(s)=C(sI-A)^{-1}B$ be a left coprime factorization
of the transfer function having the property that $\det
(sI-A)=\det D(s)$. Then the closed loop characteristic polynomial
can be written as:
\begin{equation}\label{polypol}
   \det(sI-A-BKC)=\det\vier{D(s)}{N(s)}{-K}{I}=
  \sum_\alpha g_\alpha(s)k_\alpha,
\end{equation}
where the numbers $k_\alpha$ are the Pl\"ucker coordinates (full
size minors) of the compensator $[-K\ I]$ inside
$\wedge^m\F^{m+p}$ and the polynomials $g_\alpha(s)$ are (up to
sign) the corresponding Pl\"ucker coordinates of $[D(s)\ N(s)]$.
Let ${\mathbb P}^N$ be the projective space ${\mathbb
  P}(\wedge^m\F^{m+p})$ and let
$$
E_{(A,B,C)}:=\left\{ k\in {\mathbb P}^N \mid \sum_\alpha
  g_\alpha(s)k_\alpha =0 \right\} .
$$
Since each $g_\alpha(s)$ has degree at most $n$, $E_{(A,B,C)}$
has dimension at least $N-n-1$, and its dimension equals
$N-n-1$ precisely when the $g_\alpha(s)$ span the vector space of
polynomials of degree at most $n$.  In this case, the central
projection induced by $\chi$ (see~\cite{wa96})
\begin{equation}
  \label{central}
  L_{(A,B,C)} \; :\; 
{\mathbb P}^N-E_{(A,B,C)} \  \longrightarrow \ 
  {\mathbb P}^n,\hspace{9mm}
k\  \longmapsto \ \sum_\alpha  g_\alpha(s)k_\alpha
\end{equation}
is surjective.

By\eqref{polypol}, there is a unique Pl\"ucker coordinate
$\bar{\alpha}$ with $g_{\bar{\alpha}}(s)$ of degree $n$, namely that
corresponding to the minor $\det D(s)$ of $[D(s)\ N(s)]$.
Moreover, $k_{\bar{\alpha}}=1$ and all other $g_\alpha(s)$ have
degree at most $n-1$.  Identify $\F^N\subset{\mathbb P}^N$ with
those points whose $\bar{\alpha}$th coordinate is 1.  Then the
central projection $L_{(A,B,C)}$ maps $\F^N$ to the set of monic
polynomials of degree $n$, and its complement ${\mathbb
  P}^N-\F^N$ to polynomials of degree at most $n-1$.

Every $m\times p$ compensator $K$ defines a $m$-dimensional
linear subspace of $\F^{m+p}$, the row space of $[-K\ I]$ and
therefore a point of the Grassmann variety
Grass$(m,\F^{m+p})\subset{\mathbb P}^N.$ The previous paragraph
shows this point is in $\F^N$.  Conversely, all points in
Grass$(m,\F^{m+p})\cap\F^N$ are of the form rowspace$[-K\ I]$
({\em cf.}~\cite{br81}).

The main theorem we have is:
\begin{thm}                \label{main}
  Let $\F$ be algebraically closed and $n\leq mp$.  Then the pole
  placement map $\chi_{(\bar{A},\bar{B},\bar{C})}$ is surjective
  for a particular system $(\bar{A},\bar{B},\bar{C})$ if and only
  if $\dim E_{(\bar{A},\bar{B},\bar{C})}=N-n-1$ and, for any
  $y\in \F^N-E_{(\bar{A},\bar{B},\bar{C})}\cap\F^N$,
 \begin{equation} \label{iff}
  {\rm span}\left( E_{(\bar{A},\bar{B},\bar{C})},y\right)\cap
  {\rm Grass}(m,\F^{m+p})\  \neq\  E_{(\bar{A},\bar{B},\bar{C})}\cap
  {\rm Grass}(m,\F^{m+p}).
 \end{equation}
\end{thm}
\begin{proof}
  Suppose $\chi_{(\bar{A},\bar{B},\bar{C})}$ is surjective.  Then
  the central projection $L_{(\bar{A},\bar{B},\bar{C})}$ is
  surjective and so $\dim E_{(\bar{A},\bar{B},\bar{C})}=N-n-1$.
  If for some $\hat{y}\in
  \F^N-E_{(\bar{A},\bar{B},\bar{C})}\cap\F^N,$
 \begin{equation}
   {\rm span}\left( E_{(\bar{A},\bar{B},\bar{C})},\hat{y}\right)\cap
   {\rm Grass}(m,\F^{m+p})\  = \  E_{(\bar{A},\bar{B},\bar{C})}\cap
   {\rm Grass}(m,\F^{m+p}),
 \end{equation}
 then there is also equality in\eqref{iff} for all $y\in {\rm
   span}\left( E_{(\bar{A},\bar{B},\bar{C})},\hat{y}\right).$ In
 particular, we see that the set
 $\chi_{(\bar{A},\bar{B},\bar{C})}^{-1}(L_{(\bar{A},\bar{B},
   \bar{C})}(\hat{y}))$ is empty, a contradiction.
 
 Conversely, if $\dim E_{(\bar{A},\bar{B},\bar{C})}=N-n-1$, then
 $L_{(\bar{A},\bar{B},\bar{C})}$ is surjective.  Let
 $\phi\in{\mathbb P}^n$ be any closed loop polynomial and
 $y\in{\mathbb P}^N$ satisfy
 $L_{(\bar{A},\bar{B},\bar{C})}(y)=\phi$.  Then necessarily $y\in
 \F^N$, and condition\eqref{iff} guarantees that there exists
 $P\in {\rm Grass}(m,\F^{m+p})$ with
 $L_{(\bar{A},\bar{B},\bar{C})}(P)=\phi$.  But then $P$ is the
 row space of $[-K\ I]$, for some compensator $K$.  Hence
 $\chi_{(\bar{A},\bar{B},\bar{C})}(K)=\phi$.
\end{proof}



\begin{rem}
  A system $(\bar{A},\bar{B},\bar{C})$ is {\em nondegenerate} if
  $E_{(\bar{A},\bar{B},\bar{C})}\cap {\rm
    Grass}(m,\F^{m+p})=\emptyset$. In~\cite{br81} it was shown
  that nondegenerate systems can be arbitrarily pole assigned and
  that the set of nondegenerate systems forms a generic set if
  and only if $mp\leq n$.
\end{rem}

The remainder of the paper is concerned with the question of when
the condition $mp=n$ is also sufficient for the pole placement
map $\chi$ to be generically surjective over the reals. If
$(A,B,C)$ are real matrices and if $\chi_{(A,B,C)}:
\R^{mp}\longrightarrow\R^n$ is the real pole placement map, we let
$\tilde{\chi}_{(A,B,C)}: \C^{mp}\longrightarrow\C^n$ denote the
corresponding complexified map.

\begin{thm}                  \label{criter}
  Let $\F=\R$ and assume that $mp=n$ and $d(m,p)$ is even. Then
  $\chi$ is not generically surjective if and only if there
  exists a system $(\bar{A},\bar{B},\bar{C})$ and a polynomial
  $\bar{\phi}\in\R[s]$ such that
  $\tilde{\chi}_{(\bar{A},\bar{B},\bar{C})}^{-1}(\bar{\phi})\subset\C^{mp}$
  consists of $d(m,p)$ different complex points, none of them
  real.
\end{thm}
\begin{proof}
  Assume $\chi$ is not generically surjective.  Then there exists
  a Euclidean open neighborhood $U\subset\R^{n(n+m+p)}$ for which
  $\chi_{(A,B,C)}$ is not surjective if $(A,B,C)\in U$.  Since
  $U$ is open, there exists a nondegenerate plant
  $(\bar{A},\bar{B},\bar{C})\in U$ having the property that
  $\tilde{\chi}_{(\bar{A},\bar{B},\bar{C})}^{-1}(\phi)$ consists
  of $d(m,p)$ points independent of $\phi$.  Choosing a
  polynomial $\bar{\phi}$ which is not in the image of $\chi$
  establishes one direction of the proof.
  
  On the other hand, if 
  $\tilde{\chi}_{(\bar{A},\bar{B},\bar{C})}^{-1}(\bar{\phi})\subset\C^{mp}$
  consists of $d(m,p)$ different complex points, then necessarily
  $(\bar{A},\bar{B},\bar{C})$ is a nondegenerate plant. It
  follows that there exists an open Euclidean neighborhood $U$ of
  $(\bar{A},\bar{B},\bar{C})$ consisting solely of nondegenerate systems,
  none of which can be assigned the closed loop
  characteristic polynomial $\bar{\phi}$.
\end{proof}
Theorem~\ref{criter} is interesting since it seeks a geometric
configuration where all discrete solutions are purely complex.
We use it to show that besides the case of $m=p=2$ and $n=4$,
there are other situations where $mp=n$ is not sufficient to
guarantee that $\chi$ is generically surjective over the reals.
This disproves the conjecture by S.-W. Kim mentioned in \S 1.

\begin{exmp}
  If $\F=\R$, $p=2$, $m=4$, and $n=8$ then $\chi$ is not
  generically surjective.
\end{exmp}

By Lemma 2.5, it suffices to exhibit a real system
$(\bar{A},\bar{B},\bar{C})$ and a polynomial $\bar{\phi}$ of
degree $8$ with $8$ real roots such that
$\tilde{\chi}_{(\bar{A},\bar{B},\bar{C})}^{-1}(\bar{\phi})\subset\C^8$
consists of exactly $d(4,2)=14$ purely complex solutions.  Here
is such a system:

Let $(\bar{A},\bar{B},\bar{C})$ be a minimal realization of the
system represented through a coprime factorization
$D^{-1}(s)N(s)$, where {\small
  $$
  D(s)= \left[\begin{array}{cccccc} {x}^{4}\!-\! 16
      {x}^{3}\!+\!  3 {x}^{2}\!+\! 11 x&
      \!-\! 26 {x}^{3}\!+\! 10 {x}^{2}\!+\! 7 x\!+\! 16\\
      6 {x}^{3}\!-\! 4 {x}^{2}\!-\! 9x\!-\! 5& {x}^{4}\!+\! 3
      {x}^{3}\!-\! {x}^{2}\!-\! 16 x\!-\! 13
\end{array}\right]
$$

{\scriptsize
$$
N(s)= \left[
 \begin{array}{cccccc} \!\! 9{x}^{3}\!-\! 12
    {x}^{2}\!+\! 13 x\!-\! 17& \!-\! 31 {x}^{3}\!-\! 16
    {x}^{2}\!+\! 43 x\!-\! 23& {x}^{3}\!-\! 36 {x}^{2}\!+\! 8
    x\!-\! 13&
    23 {x}^{3}\!-\! {x}^{2}\!+\! 2 x\!-\! 21\!\!\\
    8 {x}^{3}\!-\! 6{x}^{2}\!+\! 5 x\!+\! 15& 26 {x}^{3}\!-\! 14
    {x}^{2}\!-\! 11 x\!+\! 12& 11 {x}^{3}\!+\! 5 {x}^{2}\!+\! 11
    x\!+\! 33& \!-\! 7 {x}^{2}\!+\! 11x\!+\! 5
\end{array}\right]\!.
$$}
}Let
$$\bar{\phi}(s):=(s+8)(s+6)(s+4)(s+2)(s-1)(s -2)(s- 3)(s- 4).
$$

We claim that
$\tilde{\chi}_{(\bar{A},\bar{B},\bar{C})}^{-1}(\bar{\phi})$
consists of 14 purely complex solutions (displayed below).
First, we discuss how we compute
$\tilde{\chi}_{(\bar{A},\bar{B},\bar{C})}^{-1}(\bar{\phi})$ for
such a system with $n=mp$.  Identify $\C^{mp}$ with the set of
compensators $K$.  Then the $mp$ polynomial equations
\begin{equation}        \label{dets}
\det\vier{D(s)}{N(s)}{-K}{I}\ =\  0
\end{equation}
as $s$ ranges over the roots of $\bar{\phi}$ generate the ideal
of $\tilde{\chi}_{(\bar{A},\bar{B},\bar{C})}^{-1}(\bar{\phi})$ in
${\mathbb C}^{mp}$.  We used the software package
SINGULAR~\cite{SINGULAR} to compute an elimination Gr\"obner
basis~\cite{clo92} of this ideal and verify that
$\tilde{\chi}_{(\bar{A},\bar{B},\bar{C})}^{-1}(\bar{\phi})$ is
zero-dimensional with degree 14.  This calculation on the
system\eqref{dets} requires 59 seconds of CPU on a HP 9000 D250,
800 Series computer.

This Gr\"obner basis contains a univariate polynomial, the
eliminant, whose roots are the values of that variable for the
solutions.  We used the {\tt realroot} routine of Maple to
determine the number of real roots of the eliminant and {\tt
  fsolve} to compute its roots numerically.  Since we only obtain
one coordinate of each solution, we repeated this procedure to
find the others and to match the coordinates with the solutions.

Here are numerical solutions of the system\eqref{dets}.

{\tiny
  $$
\left[
\begin{array}{cc} 
    - 548.1543631072859\pm 539.02172783574002\sqrt{-1}&
      - 2966.220011381735449\pm  1301.806890926492508\sqrt{-1}\\
      227.99002317474104\mp  195.29675914098226\sqrt{-1}&
      1189.40572416765385\mp  428.190112835481936\sqrt{-1}\\
      253.619670619102274\mp  128.997418066861599\sqrt{-1}&
      1192.66093663708038\mp  127.782426659628597\sqrt{-1}\\
      -373.4608141108503\mp  376.1870941851628\sqrt{-1}&
      -907.2715490825303837\mp  2040.657619029875556\sqrt{-1}
\end{array}\right]
$$

$$
\left[
\begin{array}{cc} 
    182.1974051162797\pm  1524.2891350121054\sqrt {-1}&-
    3910.9491667600289\pm  3319.9425134666556\sqrt {-1}\\ -
    92.76689536072804\mp  494.390627883840\sqrt {-1}&
    1206.13014159582817\mp 1171.58923461208352\sqrt {-1} \\
    202.71121387564936\mp  458.78014215695346\sqrt {-1}&
    1652.30669576900037\mp  280.820983264097575\sqrt {-1}\\ -
    999.496765955436554\mp  938.918292576740638\sqrt {-1}&
    771.9810394973421\mp  4516.958140814761213\sqrt {-1}
\end{array}\right]
$$

$$
\left[
\begin{array}{cc} 
    2792.9110057318105\mp  969.00549705135278\sqrt {-1}&
    3350.9339523791667\mp  832.762320679797284\sqrt {-1}\\ -
    338.608141548768\mp  31.1420684422097\sqrt {-1}&-
    390.733153481711\mp  71.9581835765450\sqrt {-1} \\ -
    858.10666480772375\pm  463.34803831698071\sqrt {-1}&-
    1047.08493981311276\pm  448.52274247532122\sqrt {-1}\\ -
    1736.0182637110866\pm  473.54602107116131\sqrt {-1}&-
    2069.7786151738302\pm  367.88390311074763\sqrt {-1}
\end{array}\right]
$$

$$
\left[
\begin{array}{cc}
    566.14047176252718\mp  390.1690631954798\sqrt {-1}&
    894.7573009772359\mp  213.7664118348474\sqrt {-1}\\ -
    28.9144418101747\mp  8.82325220859399\sqrt {-1}&-
    31.9889032754154\mp  25.1286025912621\sqrt {-1} \\ -
    101.611268377237\pm  166.198294126534\sqrt {-1}&-
    207.075559094765\pm  158.433905818864\sqrt {-1}\\ -
    433.109410705026\pm  160.543671922194\sqrt {-1}&-
    618.358581551134\mp  8.42746099774335\sqrt {-1}
\end{array}\right]
$$

$$
\left[
\begin{array}{cc}
    - 1328.31492831596508\pm  780.43146580510958\sqrt {-1}&
    2115.8811996413627\mp  363.25099106004349\sqrt {-1}\\
    277.0599315399026\mp  134.0101686258348\sqrt {-1}&-
    426.505631447159\pm  38.4080785894925\sqrt {-1} \\
    242.753288068855\mp  128.748683783964 \sqrt {-1}&-
    380.517275415650\pm  48.1897454846160\sqrt {- 1}\\
    809.814164981704\mp  420.527784827832\sqrt {-1}&-
    1263.86094232894868\pm  149.27131835292291\sqrt {-1}
\end{array}\right ]
$$

$$
\left[
\begin{array}{cc}
  - 74.07812921055438\mp  1186.0867962658997\sqrt {-1}&
    481.83814937211068\mp  659.46539248077808\sqrt {-1}\\
    131.85311577768057\pm  223.6599712395458\sqrt {-1}&-
    28.4575338243835\pm  176.018708417247\sqrt {-1} \\
    50.0398731323218\pm  311.162560564792 \sqrt {-1}&-
    110.484321267527\pm  186.531966999705\sqrt {- 1}\\
    120.94035205524575\pm  693.23751296762126\sqrt {-1}&-
    241.138619140528\pm  419.709352592197\sqrt {-1}
\end{array}\right]
$$

$$
\left[
\begin{array}{cc}
    - 466.3420096818032\pm  2560.3776496553293\sqrt {-1}&-
    477.06216348936717\pm  1505.4573962873226\sqrt {-1}\\
    206.16217936754085\mp  504.1659905544772\sqrt {-1}&
    162.819554092696\mp  287.715806475160\sqrt {-1} \\
    198.483315335125\mp  690.317301079547 \sqrt {-1}&
    172.179197573658\mp  400.2773514799496 \sqrt {-1}\\
    350.2539156691074\mp  1658.3575908118343\sqrt {-1}&
    337.47012412920796\mp  971.424525500586678\sqrt {-1}
\end{array}\right]
$$
}

After discovering this example, we did a systematic search for
others.  In all, we generated 70 pairs $D(s),N(s)$ with random
integral polynomial entries, and, for each of the 70, considered
25 degree 8 polynomials $\bar{\phi}(s)$ with distinct integral
roots in $[-12,12]$.  Of the 1750 instances of
$\tilde{\chi}_{(\bar{A},\bar{B},\bar{C})}^{-1}(\bar{\phi})$ we
tested, {\em none} had 14 purely complex solutions, and only 3
had the `opposite' situation of 14 purely real solutions.  This
suggests that these extreme situations of real systems with real
data giving only purely complex (or purely real) solutions are
quite rare.  Despite this, we believe that it is always possible
to find such examples.  Specifically:
\begin{conj}
  If $d(m,p)$ is even and $n=mp$, then $\chi$ is not generically
  surjective over $\R$.
\end{conj}

Consider now the `opposite' situation.  Namely, for which $m,p,n$
with $n=mp$ does there exist a real system
$(\bar{A},\bar{B},\bar{C})$ and a polynomial $\bar{\phi}$ all of
whose $(n)$ roots are real such that
$\tilde{\chi}_{(\bar{A},\bar{B},\bar{C})}^{-1}(\bar{\phi})$
consists of exactly $d(m,p)$ real solutions?  Similar questions
have recently been of interest in algebraic geometry
(see~\cite{rtv95,so97dmj,so97mega} or the survey~\cite{so97sc}).
In fact, there is a precise conjecture of Shapiro and Shapiro
which is relevant to systems theory:

\begin{conj}[Shapiro-Shapiro]\label{conjSS}
  Let $(\bar{A},\bar{B},\bar{C})$ be a minimal realization of the
  system represented through a coprime factorization
  $D^{-1}(s)N(s)$, where the matrix
  $\left[D(s)\,\mid\,N(s)\right]$ has the following form: The
  first row is
  $$
  s^{m+p-1}, s^{m+p-2}, \ldots, s^2, s, 1
  $$
  and, for $1\leq j<p$, the $(j+1)$st row consists of the
  derivative of the $j$th row divided by $j$.

  Then the system is nondegenerate, and for any polynomial
  $\bar{\phi}$ of degree $mp$ with distinct real roots,
  $\tilde{\chi}_{(\bar{A},\bar{B},\bar{C})}^{-1}(\bar{\phi})$
  consists of exactly $d(m,p)$ real solutions.
\end{conj}

For example, if $m=p=3$, then we have
$$
\left[D(s)\,\mid\,N(s)\right]\ =\ \left[\begin{array}{rrr|rrr}
    s^5 & s^4 & s^3 &  s^2 & s& 1\\
    5s^4&4s^3 &3s^2 & 2s& 1& 0\\
    10s^3 &6s^2 &3s& 1&0&0\end{array}\right].
$$

For such a system, $\chi(\underline{0}) = s^{mp}$ and
$\chi^{-1}(s^{mp})=\underline{0}$, a real point with multiplicity
$d(m,p)$.  Here, $\underline{0}$ is the null compensator, the
matrix of all 0's.  Prior to learning of this conjecture, one of
us (Rosenthal) had suggested that it might be possible to perturb
$s^{mp}$ and obtain a polynomial $\bar{\phi}$ all of whose roots
are real so that $\chi^{-1}(\bar{\phi})$ consists of $d(m,p)$
real solutions.

When $p$ or $m$ is 1, this conjecture follows from
Corollary~\ref{oddR}, and when $m=p=2$, it can be verified by
hand.  All other cases remain open.  There is strong
computational evidence in support of this conjecture: In every
instance we have checked, 
$\tilde{\chi}_{(\bar{A},\bar{B},\bar{C})}^{-1}(\bar{\phi})$
consists of exactly $d(m,p)$ real solutions.  When $m=4, p=2$
(so that $d(4,2)=14$), we checked about 50 polynomials
$\bar{\phi}$.  In light of the search described above, we feel
this gives overwhelming evidence for this conjecture.  In
addition, we have considered numerous instances when $m=3,p=2$,
and a handful of instances for each of $m=5, p=2$ and $m=3, p=3$.
For each of these last two cases, $d(m,p)$ is $42$.
Unfortunately, the task of computing an elimination Gr\"obner
basis for larger $m,p$ overwhelms the HP 9000 computer we use for
these calculations.

There are other methods for solving systems of polynomials
which we have not tried, but which should work for larger $m,p$.
When $(m,p)=(5,2), (6,2)$, or $(4,3)$, we can compute a Gr\"obner
basis, and there are linear algebraic methods for solving a
polynomial system, given a Gr\"obner basis~\cite[\S 2.4]{clo97}.  Also,
homotopy continuation~\cite{ag90} algorithms which are optimized
for these systems have been developed~\cite{hss97}, and are
presently being implemented.

\begin{rem}
  The row space of the matrix $\left[D(s)\,\mid\,N(s)\right]$ of
  Conjecture~\ref{conjSS} is a $p$-plane $H(s)$ which osculates
  the moment, or rational normal curve in $\R^{m+p}$.  The
  rational normal curve is the image of the map
  $$
  s\ \longmapsto\ (s^{m+p-1}, s^{m+p-2}, \ldots, s^2, s, 1).
  $$
  
  This observation, together with the fact that all
  non-degenerate rational curves of degree $m+p-1$ in ${\mathbb
    P}^{m+p-1}$ are projectively equivalent, show that the
  conditions of Conjecture~\ref{conjSS} may be relaxed somewhat
  to the following:
  
  The row span of the matrix $\left[D(s)\,\mid\,N(s)\right]$
  equals the row span of a matrix $P(s)$ of real polynomials,
  where
\begin{enumerate}
\item The first row of $P(s)$ is a basis for all polynomials of
  degree at most $m+p-1$ and therefore defines a non-degenerate
  rational curve of degree $m+p-1$.
\item For $1\leq j<p$, the $(j+1)$st row of $P(s)$ is the
  derivative of the $j$th row of $P(s)$.
\end{enumerate}

Thus the Conjecture of Shapiro and Shapiro proposes a family of
real systems $(\bar{A},\bar{B},\bar{C})$ for which
$\tilde{\chi}_{(\bar{A},\bar{B},\bar{C})}^{-1}(\bar{\phi})$
consists of exactly $d(m,p)$ real solutions, whenever 
$\bar{\phi}$ has all real roots.

\end{rem}


\end{document}